# Power Series Method applied to Inverse Analysis in a Chemical Kinetics Problem

*E. López-Sandoval,*[*,a] *A. Mello,*[b] *J. J. Godina- Nava*[c] *and A. R. Samana*[a]

[a]*Departamento de Ciências Exatas e Tecnológicas, Universidade Estadual de Santa Cruz,*

*Campus Soane Nazaré de Andrade, Rodovia Jorge Amado, Km 16, Bairro Salobrinho,*

*45662-900 Ilhéus, BA, Brazil*

[b] *Centro Brasileiro de Pesquisas Físicas,*

*Rua Dr. Xavier Sigaud, 150*

*CEP 22290-180, Rio de Janeiro, RJ, Brazil.*

[c]*Departamento de Física,*

*Centro de Investigación y de Estudios Avanzados*

*del Instituto Politecnico Nacional*

*Ap. Postal 14-740, 07000, México, D. F. México*

*e-mail: sandoval@cbpf.br

**Abstract**

We present a deterministic and direct method based on Power Series Method applied to solve a chemical kinetics inverse problem. We use a power series as a trial solution approximation, where the data are used as initial condition and boundary value, to solve a coupled differential equation of chemical kinetics obtaining the rate constant parameters.

**Keywords**: Chemical kinetics, Rate Parameters, Inverse Problem, Power Series Method.

## 1. Introduction

A problem generally is ill-posed if it does not comply one of three conditions: (i) existence, (ii) uniqueness and, (iii) continuity.[1] Inverse analysis problems (IAP) are ill-posed because we do not know the physical parameter of the model and, therefore we cannot solve its differential equations in a direct way. It means that there is not a unique solution because we have more variables than equations. Therefore, it is necessary to employ stochastic methods, using trial values of the unknown parameters as input values in the Differential Equations (DE) of the model, solving via some numerical method until the difference between experimental data and numerical results is minimized. Such a method is an indirect or inverse one, where the difficulty is in their non linear nature and irregularity. Stochastic methods such as Particle Swarm Optimization,[2-4] Genetic algorithm[4] or Neural Network[5,6] are robust and efficient approaches to find the global optimum. Nevertheless, they have problems with many parameters because of the high computational cost.

Another method is a deterministic one, known as Non Linear Regression (NLR) and it consists in the optimization of an objective function, the mean square error.[4, 5] To optimize this function, exist many methods, but the best known are the method of steepest descent, the method of Gauss-Newton and the Marquardt method, that is a hybrid method between two previous. [4, 5] These methods are efficient, but they have a difficulty to find the global optimum when exist very much parameters, and therefore it demands a high computational cost. As M. C. Colaço et al says:[4] Deterministic methods are in general computationally faster than stochastic methods, although they can converge to a local minima or maxima, instead of the global one. On the other hand, stochastic algorithms can ideally converge to a global maxima or minima, although they are computationally slower than the deterministic ones.

## 2. Methodology

Firstly, we will assume a polynomial regression of the experimental datum to avoid include experimental errors. With these quasi-exact datum we will apply an interpolation, evaluating each one of the experimental datum of the chemical kinetics system in the trial polynomial solution. Besides, we have to evaluate this polynomial in the differential equation (DE) for each data to obtain its curvature (so we obtain a better fitting because we are using both restrictions). In this way, we use the Power Series Method (PSM), as a solution of the differential equation of our chemical kinetic model. Solving an algebraic system equation with the same number that variables, we will obtain the chemical parameters. Our method was inspired of the Spectral Collocation Method.[6-9]

The Power Series Method (PSM) is a technique which has been often used to solve linear ordinary differential equations (ODE),[10, 11] and partial differential equations (PDE).[12, 13] However, as we know, the PSM is also actually a powerful technique to solve non linear PDE.[14-18] We use the data on the concentration of the substance obtained experimentally, as the initial conditions and boundary values.

The PSM is a proposal to find a semi-analytic solution as an asymptotic approximation (in space and time) of a finite series with a minimal error in the expansion terms. We know that almost the totality of the non linear PDE (NLPDE) have not a solution with an analytic expression. The proposal is to construct a solution using the power series, taking advantage of the capacity of this method to represent any function with a polynomial, building an approximate solution:

$$f(x_1,...,x_d) = \sum_{n_1=0}^{\infty}\sum_{n_2=0}^{\infty}...\sum_{n_d}^{\infty} a_{n_1...n_d} (x_1 - a_1)^{n_1}...(x_d - a_d)^{n_d} . \qquad (1)$$

In a similar way than the PSM[12-18], this possible solution is constructed with this class of trial polynomial function (1).

We know that the chemical kinetic problems have a DE similar to:

$$dX_i(t)/dt = F_i(k_j, X_i(t)),  \qquad (2)$$

where $i=1,2,...n$ are the indices of the coupled linear differential equation system, with $j=1,2,...,m$ are the indices of the unknown parameters. According to the method, we can adopt a function solution as a power series function:

$$X_i(t) \approx \sum_{j=0}^{n} a_j t^j . \qquad (3)$$

Substituting (3) into the DE (2),

$$\sum_{j=1}^{n} j a_j t^{j-1} = F_i(k_j, \sum_{j=0}^{n} a_j t^j), \qquad (4)$$

we obtain an equation system with $n$ equations and $n+m$ unknown variables. Thus, we have an ill-posed problem because we have more variables than equations. To achieve a well posed problem[1], we have to include the data of the concentration evaluated in the trial solution (3) for each time until that $m$ additional equations are obtained.

The polynomial approximation must be of a certain degree yielding the best fitting of the experimental data concentrations. However, the complete number of restrictions with the trial function evaluation help us to perform a better fitting and to find the physical parameters more exactly.

### 3. Results and Discussion

A chemical kinetic inverse problem of three-species reversible isomerization is proposed, with six rate constants parameters of the product concentrations[19]:

$$\frac{dC_a}{dt} = -k_{12}C_a + k_{21}C_b - k_{13}C_a + k_{31}C_c,$$

$$\frac{dC_b}{dt} = -k_{21}C_b + k_{12}C_a - k_{23}C_b + k_{32}C_c, \qquad (5)$$

$$\frac{dC_c}{dt} = -k_{31}C_c + k_{13}C_a - k_{32}C_c + k_{23}C_b.$$

We want to obtain the rate constant parameters $k_{12}$, $k_{21}$, $k_{13}$, $k_{31}$, $k_{32}$, $k_{23}$, knowing the values of the initial concentrations, and the concentrations of the three substances $C_a$ $C_b$ and $C_c$, for each time of the reactions, as is shown in Table 1, obtained from Reference. [19]

**Table 1.** Time evolution concentrations of the three substances $C_a$, $C_b$ and $C_c$ (we presented rounded values).

| t/ Concentration | $C_a$ | $C_b$ | $C_c$ |
|---|---|---|---|
| 0 | 0.04 | 0.96 | 0 |
| 0.05 | 0.3590 | 0.5373 | 0.1037 |
| 0.1 | 0.4886 | 0.3703 | 0.1410 |
| 0.15 | 0.5413 | 0.3053 | 0.1534 |

The next power series are proposed as trials solutions of the differential equation system (5):

$$C_a(t) = \sum_{n=0}^{10} a_n t^n, \qquad C_b(t) = \sum_{n=0}^{10} b_n t^n, \qquad C_c(t) = \sum_{n=0}^{10} c_n t^n. \qquad (6)$$

and their respective derivatives:

$$\frac{dC_a}{dt} = \sum_{n=0}^{10} n a_n t^{n-1}, \qquad \frac{dC_b}{dt} = \sum_{n=0}^{10} n b_n t^{n-1}, \qquad \frac{dC_c}{dt} = \sum_{n=0}^{10} n c_n t^{n-1}. \qquad (7)$$

Substituting these (6) and (7) into the differential equations (5), we obtain the recurrence relation:

$$a_{n+1} = \frac{1}{(n+1)}\left[-a_n(k_{12}+k_{13}) + b_n k_{21} + c_n k_{31}\right],$$

$$b_{n+1} = \frac{1}{(n+1)}\left[-b_n(k_{21}+k_{23}) + a_n k_{12} + c_n k_{32}\right], \qquad (8)$$

$$c_{n+1} = \frac{1}{(n+1)}\left[-c_n(k_{31}+k_{32}) + a_n k_{13} + b_n k_{23}\right].$$

We obtain a system of equations when $n=0$ until $n=10$. Here, we have 30 equations with 33 unknown parameters ($a_i$, $b_i$ and $c_i$, with $i=0, 1,…10$), and 6 chemical parameters ($k_{12}$, $k_{21}$, $k_{31}$, $k_{13}$, $k_{32}$, $k_{23}$). So there, we need more 9 additional equations. We can to obtain these ones with the data in Table 2 substituting in the trial solutions (6):

For t=0:

$$C_a(t=0) = \left(\sum_{n=0}^{10} a_n t^n\right)_{t=0} = a_0 = 0.04,$$

$$C_b(t=0) = \left(\sum_{n=0}^{10} b_n t^n\right)_{t=0} = b_0 = 0.96,$$

$$C_c(t=0) = \left(\sum_{n=0}^{10} c_n t^n\right)_{t=0} = c_0 = 0.$$

For t=0.05:

$$C_a(t=0.05) = \left(\sum_{n=0}^{10} a_n t^n\right)_{t=0.05} = 0.3590,$$

$$C_b(t=0.05) = \left(\sum_{n=0}^{10} b_n t^n\right)_{t=0.05} = 0.5373, \qquad (9)$$

For t=0.10:

$$C_a(t = 0.10) = \left(\sum_{n=0}^{10} a_n t^n\right)_{t=0.10} = 0.4886,$$

$$C_b(t = 0.10) = \left(\sum_{n=0}^{10} b_n t^n\right)_{t=0.10} = 0.3703.$$

For t=0.15:

$$C_a(t = 0.15) = \left(\sum_{n=0}^{10} a_n t^n\right)_{t=0.15} = 0.5413,$$

$$C_b(t = 0.15) = \left(\sum_{n=0}^{10} b_n t^n\right)_{t=0.15} = 0.3053.$$

Until *t=0.15* it was enough to obtain a system of equation with the same number of unknown variables. So then, we take only 9 values of the concentration because only we need 9 equations for to obtain a well posed problem in the determination of the 6 rate constants parameters. Also, we use only 2 functions because it was necessary to avoid redundancy in the solution. We solve such a system of non linear algebraic equations, using the Marquardt method with Matlab® software. Thus, we obtain the parameters of the function solution as it shown in the table 2. Fig. 1 shows the graphical representation of this function. The calculated rate constants parameter values are presented in Table 3.

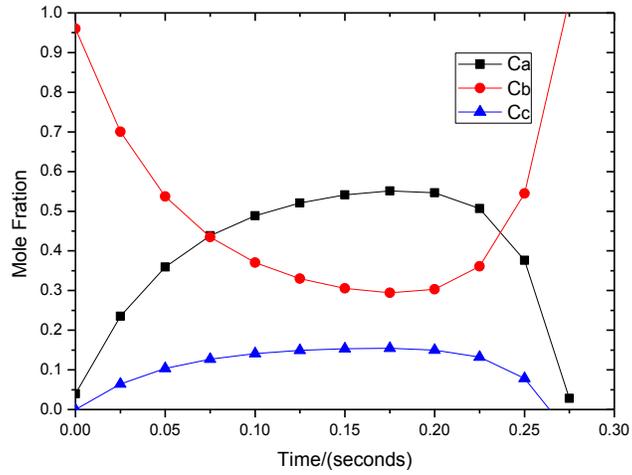

**Figure 1.** Trajectories of state variables for the reaction of three substances.

**Table 2.** Parameters of function solution for the reaction of three substances (we presented rounded values).

| $i$ | $a_i$ | $b_i$ | $c_i$ |
|---|---|---|---|
| 1 | 9.7053 | -12.98148 | 3.2761 |
| 2 | -88.3265 | 120.6538 | -32.3273 |
| 3 | 541.4386 | -748.2168 | 206.7781 |
| 4 | -2503.5057 | 3481.5326 | -978.02696 |
| 5 | 9289.6028 | -12963.10536 | 3673.5025 |
| 6 | -28774.2791 | 40227.5356 | -11453.2565 |
| 7 | 76465.6390 | -107009.376 | 30543.7370 |
| 8 | -177890.799 | 249083.4838 | -71192.6840 |
| 9 | 367964.2187 | -515376.3624 | 147412.1436 |
| 10 | -685115.307 | 959735.1993 | -274619.8921 |

**Table 3.** Rate constant parameters obtained of the reaction of three substances.

|  | $k_{12}$ | $k_{21}$ | $k_{23}$ | $k_{32}$ | $k_{31}$ | $k_{13}$ |
|---|---|---|---|---|---|---|
| Exact | 4.623 | 10.344 | 3.371 | 5.616 | 3.724 | 1.0 |
| Calculated | 4.62299... | 10.34399... | 3.3709... | 5.6160... | 3.72399... | 1.000... |
| Absolute Error | 1.28e-011 | 1.247e-013 | 1.825e-013 | 4.31e-011 | 1.161e-011 | 9.273e-012 |

In a similar way, we will solve a four-species reversible isomerization problem. Here, we want to calculate the eight rate constants parameters of the DE of four concentrations of different substances[19]:

$$\frac{dC_a}{dt} = -k_{12}C_a + k_{21}C_b - k_{14}C_a + k_{41}C_d,$$

$$\frac{dC_b}{dt} = -k_{21}C_b + k_{12}C_a - k_{23}C_b + k_{32}C_c, \quad (10)$$

$$\frac{dC_c}{dt} = -k_{34}C_c + k_{43}C_d - k_{32}C_c + k_{23}C_b,$$

$$\frac{dC_d}{dt} = -k_{43}C_d + k_{34}C_c - k_{41}C_d + k_{41}C_a,$$

and to find the rate constants parameters $k_{12}$, $k_{21}$, $k_{14}$, $k_{41}$, $k_{32}$, $k_{23}$, $k_{34}$, $k_{43}$. The values of the concentrations of the four substances, $C_a(t)$, $C_b(t)$, $C_c(t)$ and $C_d(t)$ for three different times of its concentration are shown in the Table 4, obtained from Reference.[19]

**Table 4.** Time evolution concentrations of the four substances $C_a$ $C_b$, $C_c$ and $C_d$. (we presented rounded values).

| t/Concentration | $C_a$ | $C_b$ | $C_c$ | $C_d$ |
|---|---|---|---|---|
| 0 | 0.4 | 0.2 | 0.1 | 0.3 |
| 0.025 | 0.3979 | 0.1980 | 0.1206 | 0.2835 |
| 0.05 | 0.3931 | 0.2002 | 0.1348 | 0.2719 |
| 0.075 | 0.3873 | 0.2045 | 0.1448 | 0.2633 |

The next power series are proposed as trials solutions:

$$C_a(t) = \sum_{n=0}^{8} a_n t^n, \quad C_b(t) = \sum_{n=0}^{8} b_n t^n, \quad C_c(t) = \sum_{n=0}^{8} c_n t^n, \quad C_d(t) = \sum_{n=0}^{10} d_n t^n \quad (11)$$

and their respective derivatives:

$$\frac{dC_a}{dt} = \sum_{n=0}^{10} na_n t^{n-1}, \quad \frac{dC_b}{dt} = \sum_{n=0}^{10} nb_n t^{n-1}, \quad \frac{dC_c}{dt} = \sum_{n=0}^{10} nc_n t^{n-1}, \quad \frac{dC_d}{dt} = \sum_{n=0}^{10} nd_n t^{n-1} \quad (12)$$

Substituting the trial solutions, (8) and (9), into the differential equations (17), we obtain the recurrence relation:

$$a_{n+1} = \frac{1}{(n+1)}\left[-a_n(k_{12} + k_{14}) + b_n k_{21} + d_n k_{41}\right],$$

$$b_{n+1} = \frac{1}{(n+1)}\left[-b_n(k_{21} + k_{23}) + a_n k_{12} + c_n k_{32}\right],$$

$$c_{n+1} = \frac{1}{(n+1)}\left[-c_n(k_{32} + k_{34}) + d_n k_{43} + b_n k_{23}\right], \quad (13)$$

$$d_{n+1} = \frac{1}{(n+1)}\left[-d_n(k_{41} + k_{43}) + c_n k_{34} + a_n k_{41}\right],$$

when $n=0$ until $n=8$, we obtain the system of equations. Here, we have 32 equations with 36 unknown parameters ($a_i$, $b_i$, $c_i$ and $d_i$, with $i=0, 1,...,8$) and 8 unknown chemical parameters ($k_{12}$, $k_{21}$, $k_{14}$, $k_{41}$, $k_{23}$, $k_{32}$, $k_{34}$, $k_{43}$). So there, we need more 12 additional equations. First, we have the data of the Table 4 that are considered as the boundary values. Evaluating the trial series solution (11) with these data:

For t=0:

$$C_a(t=0) = \left(\sum_{n=0}^{8} a_n t^n\right)_{t=0} = a_0 = 0.4,$$

$$C_b(t=0) = \left(\sum_{n=0}^{8} b_n t^n\right)_{t=0} = b_0 = 0.2,$$

$$C_c(t=0) = \left(\sum_{n=0}^{8} c_n t^n\right)_{t=0} = c_0 = 0.1,$$

$$C_d(t=0) = \left(\sum_{n=0}^{8} d_n t^n\right)_{t=0} = d_0 = 0.3.$$

For t=0.025:

$$C_a(t=0.025) = \left(\sum_{n=0}^{8} a_n t^n\right)_{t=0.025} = 0.3979,$$

$$C_b(t=0.025) = \left(\sum_{n=0}^{8} b_n t^n\right)_{t=0.025} = 0.1980,$$

$$C_c(t=0.025) = \left(\sum_{n=0}^{8} c_n t^n\right)_{t=0.025} = 0.1206,$$

For t=0.05:

$$C_a(t=0.05) = \left(\sum_{n=0}^{8} a_n t^n\right)_{t=0.05} = 0.3931,$$

$$C_b(t=0.05) = \left(\sum_{n=0}^{8} b_n t^n\right)_{t=0.05} = 0.2002,$$

$$C_c(t=0.05) = \left(\sum_{n=0}^{8} c_n t^n\right)_{t=0.05} = 0.1348.$$

For t=0.075:

$$C_a(t=0.075) = \left(\sum_{n=0}^{8} a_n t^n\right)_{t=0.075} = 0.3873,$$

$$C_b(t=0.075) = \left(\sum_{n=0}^{8} b_n t^n\right)_{t=0.075} = 0.2045.$$

Until t=0.075 it was enough to obtain a system of equation with the same number of unknown variables. We take only 12 values of the concentration, because only we need 12 equations to obtain a well posed problem in the determination of the 8 rate constants parameters. Also, we use only 3 functions because it was necessary to avoid redundancy in the solution.

Again, solving with Matlab®, we obtain the coefficients of the functions of the substance concentration as it is shown in Table 5. Fig. 2 shows the graphical representation of this function. The chemical constant parameter values are shown in Table 6.

**Table 5.** Parameters of function solution for the reaction of four substances (we presented rounded values).

| $i$ | $a_i$ | $b_i$ | $c_i$ | $d_i$ |
|---|---|---|---|---|
| 1 | 1.5056e-4 | 0.2001 | 0.99996 | -0.8800 |
| 2 | -4.3107 | 5.8057 | -8.0977 | 7.2629 |
| 3 | 40.9269 | -44.0145 | 50.8307 | -52.5173 |
| 4 | -245.5467 | 225.9489 | -255.0442 | 302.1062 |
| 5 | 1131.6350 | -923.3400 | 1056.9805 | -1391.8031 |
| 6 | -4285.2524 | 3190.6599 | -3724.4624 | 5300.9604 |
| 7 | 13815.2987 | -9607.7826 | 11410.8308 | -17180.1817 |
| 8 | -38815.5232 | 25675.2463 | -30913.6774 | 48459.3498 |

**Table 6.** Rate constant parameters obtained of the reaction of four substances.

| | $k_{12}$ | $k_{21}$ | $k_{23}$ | $k_{32}$ | $k_{34}$ | $k_{43}$ | $k_{41}$ | $k_{14}$ |
|---|---|---|---|---|---|---|---|---|
| *Exact* | *1* | *3* | *5* | *10* | *2* | *4* | *10* | *8* |
| *Calculated Parameter* | *1.01659* | *3.0553* | *4.98027* | *10.00394* | *1.98722* | *4.01008* | *10.01045* | *8.01851* |
| *Absolute error* | *0.01659* | *0.0553* | *0.01972* | *0.00394* | *0.01277* | *0.01008* | *0.01045* | *0.01851* |

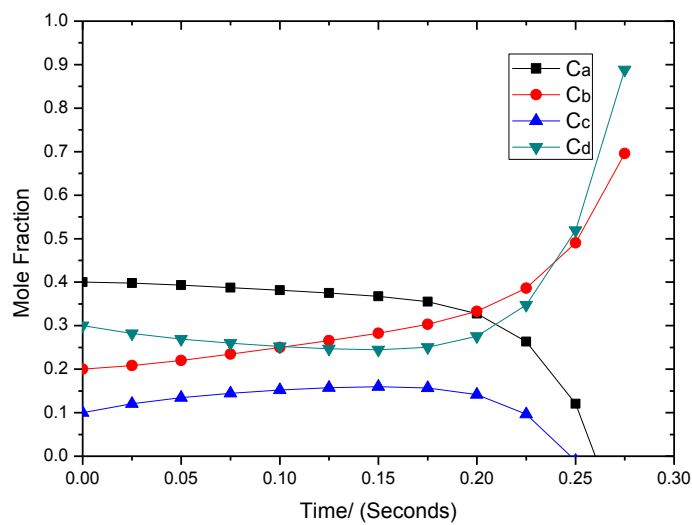

**Figure 2.** Trajectories of state variables for the reaction of four substances.

In this letter, we have shown that is possible to solve a chemical kinetics problem in a direct and deterministic way using the Power Series Method. With this method, it is possible to obtain a well posed problem evaluating the data of the concentrations obtained from its graphical representation in the polynomial series and in the DE, resulting in the same number of equations as well as variables. With these restrictions, we assure that the trial solution is well posed, complying the conditions of: (i) existence, (ii) uniqueness and, (iii) continuity. [1]

In NLR there are many optimization deterministic methods, then an exhausted search is necessary in the phase space that depends of the number of parameters. Then it is very difficult to find the global minimum because the search space is huge, and it could have very much local minima. Then, it could converge slowly, or not converge, or converges in an inexact way because, as Motulsky et al[5] say, the function becomes ill-conditioned or the Jacobian matrix becomes singular. This occurs in several situations: 1) The data contain numbers that are too large or too small. 2) The selected equation does not reasonably fit the data. 3) The initial values are far from correct. 4) The data points are quite scattered. 5) The data were not collected over a sufficiently wide range of x values. 6) The computer calculations were not sufficiently precise (not enough significant digits).

With the PSM, we avoid some of these problems because the solution does a very good fitting, where we used all the necessary constrains for a better curve fitting. On the other hand, the PSM is very precise to solve DE and the first and second derivative are obtained more easy and more exact than with other numerical method.[18]

## 4. Conclusion

The PSM is a semi-analytic technique which obtains, in an easier way, the solution of difficult differential equations with an approximated closed form expression.[18] Furthermore, the PSM is especially useful to solve non-linear equations, opening the possibility to treat a more complex problem of inverse analysis. The solution could be developed until the necessary degree into the power series that could be to fit to the data, obtaining an approximated value of the parameters.

**Acknowledgements**

E.L.S. et al. are fellows to Roman Lopez Sandoval and Rogelio Ospina for their very careful reading and helpful discussion of the manuscript. We also fellows the support of Brazilian agency CNPq 680.023/2008-9 (PCI DTI-7A) and 170.191/2011-7 (PCI-BEV). J.J. Godina-Nava fellows CONACYT-México 132432 and 254931. E.L.S and A.R.S fellows the support of brazilian agency PNPD-CAPES 201330748-28007018010P2.